\documentclass[10pt]{amsart}

 \usepackage{amssymb,latexsym,amsmath,amsthm}
\usepackage{graphicx} 

\usepackage[numbers,sort&compress]{natbib}

\usepackage{fullpage}

 \renewcommand{\div}{\mathop{\mathrm{div}}\nolimits}
 \usepackage[numbers,sort&compress]{natbib}

\newtheorem*{thm*}{Theorem A}

\newtheorem{thm}{Theorem}[section]

\newtheorem{notation}{Notation}[section]

\newtheorem{lemma}{Lemma}[section]

\newtheorem{cor}{Corollary}[section]

\numberwithin{equation}{section}

\usepackage{graphicx} 

\usepackage{epsfig}
\usepackage{graphicx}
\usepackage{caption}
\usepackage{subcaption}
\usepackage{float}
\usepackage{amsthm}
\usepackage{amssymb}
\usepackage{booktabs} 
\usepackage{array} 
\usepackage{paralist} 
\usepackage{verbatim} 
\usepackage{subfig} 
\usepackage{fancyhdr} 
\begin{document}

\author{Mostafa Fazly}

\address{Department of Mathematics, The University of Texas at San Antonio, San Antonio, TX 78249, USA}
\email{mostafa.fazly@utsa.edu}

\author{Mark Lewis}
\address{Department of Mathematical  \& Statistical Sciences,  University of Alberta, Edmonton, AB T6G 2G1 Canada.}
\email{mark.lewis@ualberta.ca}

\author{Hao Wang}

\address{Department of Mathematical  \& Statistical Sciences,  University of Alberta, Edmonton, AB T6G 2G1 Canada.}
\email{hao8@ualberta.ca}

\thanks{MAL gratefully acknowledges a Canada Research Chair and NSERC Discovery grants. MF and HW gratefully acknowledge NSERC Discovery grants.}

\def\IR{{\mathbb{R}}}

\newcommand{ \un }{\underline}

\newcommand{ \bg }{\begin{equation}}

\newcommand{ \ed }{\end{equation}}

\newcommand{ \Om }{ \Omega}

\newcommand{ \pOm}{\partial \Omega}

\title{On impulsive reaction-diffusion models in higher dimensions}
\maketitle

\begin{abstract}
Assume that  $N_m(x)$ denotes the density of the population at a point $x$ at the beginning of the reproductive season in the $m$th year. We study the following impulsive reaction-diffusion model for any $m\in \mathbb Z^+$
\begin{eqnarray*}\label{}
 \ \ \ \ \  \left\{ \begin{array}{lcl}
    u^{(m)}_t  = \div(A\nabla u^{(m)}-a u^{(m)})   + f(u^{(m)})  \quad  \text{for} \ \  (x,t)\in\Omega\times (0,1] \\
 u^{(m)}(x,0)=g(N_m(x))  \quad  \text{for} \ \ x\in \Omega\\
  N_{m+1}(x):=u^{(m)}(x,1)   \quad  \text{for} \ \ x\in \Omega
\end{array}\right.
  \end{eqnarray*}
for  functions $f,g$,  a drift $a$ and  a diffusion matrix $A$ and $\Omega\subset \mathbb R^n$.  Study of this model requires a simultaneous analysis of the differential equation and the recurrence relation.    When boundary conditions are hostile we provide critical domain results showing how extinction versus persistence of the species arises,  depending on the size and geometry of the domain.
We show that there exists an {\it extreme volume size} such that if $|\Omega|$ falls below this size  the species is driven extinct,    regardless of the geometry of the domain.  To construct such  extreme volume sizes and  critical domain sizes,  we apply Schwarz symmetrization rearrangement arguments,  the classical Rayleigh-Faber-Krahn inequality  and  the spectrum of uniformly elliptic operators.    The critical domain results provide qualitative insight regarding long-term dynamics for the model.  Lastly, we provide applications of our main results to certain biological reaction-diffusion models regarding  marine reserve, terrestrial reserve, insect pest outbreak and  population subject to climate change. 
\end{abstract}

\noindent
{\it \footnotesize 2010 Mathematics Subject Classification}. {\scriptsize  92B05, 35K57, 92D40, 92D50, 92D25.}\\
{\it \footnotesize Keywords: Reaction-diffusion models, persistence versus extinction of species, eigenvalue problems for differential operators, rearrangement arguments, population dynamics}. {\scriptsize }


\section{Introduction}
Impulsive reaction-diffusion equation models for species with distinct reproductive
and dispersal stages were proposed by Lewis and Li in \cite{LL}. These models can be considered to be description for a seasonal birth pulse
plus nonlinear mortality and dispersal throughout the year. Alternatively,   they
can describe seasonal harvesting, plus nonlinear birth and mortality as well
as dispersal throughout the year.  The population of a species at the beginning of year $m$ is denoted by  $N_m(x)$. We assume that reproduction (or harvesting) occurs at the beginning of the year, via a discrete time map, $g$, after which there is  birth, mortality and dispersal via a reaction-diffusion  for a population with density $u^{(m)}(x,t)$. At the end of this year the density $u^{(m)}(x,1)$ provides the population density for the start of year $m+1$, $N_{m+1}(x)$.   We examine solutions of the following system for any $m\in \mathbb Z^+$
\begin{eqnarray}\label{main}
   \ \ \ \  \left\{ \begin{array}{lcl}
  u^{(m)}_t = \div(A\nabla u^{(m)}-a u^{(m)})  + f(u^{(m)})  \quad  \text{for} \ \  (x,t)\in\Omega\times (0,1]   , \\
u^{(m)}(x,t)=0   \quad  \text{for} \ \  (x,t)\in\partial\Omega\times (0,1]    ,  \\
 u^{(m)}(x,0)=g(N_m(x))  \quad  \text{for} \ \ x\in \Omega   , \\
  N_{m+1}(x):=u^{(m)}(x,1)   \quad  \text{for} \ \ x\in \Omega,
\end{array}\right.
  \end{eqnarray}
where $\Omega$ is a set in $\mathbb{R}^n$,  $A$ is a constant symmetric positive definite matrix and $ a$ is a constant vector.  We assume that the function $g$ satisfies the following assumptions;
 \begin{enumerate}
 \item[(G0)] {\it $g$ is a continuous positive function in $\mathbb R^+$ and $g(0)=0< g'(0)$ and there exists a positive constant $M\in (0,\infty]$ such that  $g(N)$ is nondecreasing for $0<N \le M$. In other words, monotonicity of $g$ is required only on a subset of $\mathbb R^+$.}
 \item[(G1)] {\it  there exists a positive constant $\bar M \le M $ such that $g(N)\le g'(0)N$ for $0<N<\bar M$.} 
 \item[(G2)]{\it  there exists a differentiable function $h$ for $h(0)=h'(0)=0$ and a constant $\tilde M\le M$ so that 
 $g(N)\ge g'(0)N-h(N)$ for $0<N< \tilde M$. }
 \end{enumerate}
   Note that  the linear function 
\begin{equation}\label{lg}
g(N)=b N, 
\end{equation}
 where $b$ is a positive constant satisfies all of the above assumptions.  For the above $g$ when $b=1$ and $f(N)=N(1-N)$,  model (\ref{main}) recovers the classical Fisher's equation introduced  by Fisher in \cite{fish} and Kolmogorov, Petrowsky, and Piscounoff (KPP) in \cite{kpp} in 1937.    One can consider nonlinear functions for $g$ such as the Ricker function,  that is 
\begin{equation}\label{rg}
g(N)=N e^{r(1-N)}, 
\end{equation} where $r$ is a positive constant.  For the optimal stocking rates for fisheries  mathematical biologists often apply  the Ricker model \cite{rick} introduced in 1954  to study  salmon populations with scramble competition for spawning sites leading to overcompensatory dynamics. 
The Ricker function is nondecreasing for $0<N \le M=\frac{1}{r}$ and satisfies all of the assumptions.     Note also that the  Beverton-Holt function 
\begin{equation}\label{bg}
g(N) =\frac{(1+\lambda)N}{1+\lambda N}, 
\end{equation}  for positive constant $\lambda$  is an increasing function.    
 The Beverton-Holt model was introduced to understand the dynamics of compensatory competition in  
 fisheries by Beverton and Holt \cite{bh} in 1957.     Another example is the Skellam function 
\begin{equation}\label{sg}
g(N)=R(1-e^{-bN}), 
\end{equation}   where $R $ and $b$ are positive constants. This function was introduced by Skellam in 1951 in \cite{s} to study population density for animals, such as birds, which have contest competition for nesting sites that leads to compensatory competition dynamics. 
Note that the Skellam function behaves similar to the  Beverton-Holt function and it is nondecreasing for any $N>0$.  We shall use these functions in the application section (Section \ref{app}).  We refer interested readers to \cite{st13} for more functional forms with biological applications.   We now provide some assumptions on the function $f$. We suppose that 
 \begin{enumerate}
 \item[(F0)]{\it  $f(.)$ is a continuous function and  $f(0)=0$. We also assume that $f'(0)\neq 0$.}
 \item[(F1)] {\it  there exists a differentiable function $h$ for $h(0)=h'(0)=0$  so that 
 $$ f'(0)N-h(N) \le f(N)\le f'(0) N \ \ \ \text{for}\ \ \ N \in \mathbb R^+.$$ }
 \end{enumerate}
 Note that we do not have any assumption on the sign of $f(.)$ and $f'(0)$.    Note that $f(N)=N(1-N)$,  $f(N)=bN$ for $b\in\mathbb R$ and $f(N)=\alpha N - \beta N^2$ for $\alpha \in \mathbb R$ and $\beta\in\mathbb R^+$ satisfy the above assumptions (F0) and (F1).

  Suppose that $a=0$ and $A=0$ then $u(x,t)$ only depends on time and not space, meaning that individuals do not advect or diffuse. Assume that $N_m$ represents the
number of individuals at the beginning of reproductive stage in the $m$th year. Then
 \begin{eqnarray}\label{ut}
 \left\{ \begin{array}{lcl}
    u_t (t)= f(u(t))  \quad  \text{for} \ \  t\in (0,1] , \\
 u(0)=g(N_m)    , \\
  N_{m+1}:=u(1).
\end{array}\right.
  \end{eqnarray}
   Separation of variables shows that
   \begin{equation}\label{equiseq}
  \int_{g(N_m)}^{N_{m+1}} \frac{d \omega}{f(\omega)}=1.
  \end{equation}
Note that a positive constant equilibrium of (\ref{ut})  satisfies
 \begin{equation}\label{equib}
 \int_{g(N)}^{N} \frac{d \omega}{f(\omega)}=1.
\end{equation}
Assume that $f$ satisfies (F0)-(F1) and $g$ satisfies (G0)-(G2) then
 \begin{equation}\label{Nstarine}
 1=\int_{g(N)}^{N} \frac{d \omega}{f(\omega)}\ge  \frac{1}{f'(0)}  \int_{g(N)}^{N} \frac{d \omega}{\omega}=  \frac{1}{f'(0)}  \ln \left |\frac{N}{g(N)} \right| \ge \frac{1}{f'(0)}  \ln \left |\frac{1}{g'(0)} \right|  .
   \end{equation}
In the light of above computations,   we assume that   
 \begin{equation}\label{conditionfg}
 e^{f'(0)} g'(0)>1,
\end{equation}
 and an $N^*>0$ exists such that $f\neq 0$ on the closed interval with endpoints $N^*$ and $g(N^*)$ and 
 \begin{equation}\label{nstar} 
 \int_{g(N^*)}^{N^*} \frac{d \omega}{f(\omega)}=1.
\end{equation}
We also assume that $g$ is nondecreasing on $[0,N^*]$; that is equivalent to considering $M:=N^*$ in the assumption (G0).  Note that for the case of equality in (\ref{Nstarine}) there might only be the zero equilibrium,   for examples, see (\ref{Nstar}) and (\ref{Nstarr}) in the application section (Section \ref{app}).

 The above equation (\ref{main}) defines a recurrence relation for $N_m(x)$ as
  \begin{equation}\label{op}
  N_{m+1}(x)=Q[N_m(x)]   \quad  \text{for} \ \ x\in \Omega\subset\mathbb{R}^n,
  \end{equation}
  where $m\ge 0$ and $Q$ is an operator that depends on $A, a,f,g$.  While most of the results provided in this paper are valid in any dimensions, we shall focus on the case of $n\le 3$ for applications. For notational convenience we drop superscript $(m)$ for $u^{(m)}(x,t)$, rewriting it as $u(x,t)$. The remarkable point about the impulsive reaction-diffusion equation (\ref{main})-(\ref{op}) is that it is a mixture of a differential equation and a recurrence relation.   Therefore, one may expect that the analysis of this model requires a simultaneous analysis of the continuous  and  discrete  type. If the yearlong activities are modelled by impulsive dynamical systems,  models of the form (\ref{main}) have a longer history and have been given various names: discrete time metered models \cite{cl}, sequential-continuous models \cite{bc} and semi-discrete models \cite{sn}. We refer interested readers to \cite{st11} for more information.  
  
 In this paper, we provide critical domain size for extinction versus persistence of populations for impulsive reaction-diffusion models of the form (\ref{main}) defined on domain $\Omega\subset\mathbb R^n$.  It is known that the geometry of the domain has fundamental impacts on qualitative behaviour of solutions of equations in higher dimensions.     We consider domains with various geometric structures in dimensions $n\ge1$ including convex and concave domains and also domains with smooth and non smooth boundaries. 
 
 The discrete time models of the form (\ref{op}) are studied extensively in the literature, in particular in the foundational work of Weinberger \cite{w82}, where $N_m(x)$ represents the gene fraction or population density at time $n$ at the point $x$ of the habitat and $Q$ is an operator  on a certain set of functions on the habitat. It is shown by Weinberger \cite{w82} that under a few biologically reasonable hypothesis on the operator $Q$ the results similar to those for  the Fisher and Kolmogoroff, Petrowsky, and Piscounoff (KPP) types for models (\ref{op}) hold. In other words, given a direction vector $e$,  the recurrence relation (\ref{op}) admits a nonincreasing planar traveling wave solution   for every $c\ge c^*(e)$ and, more importantly,  there will be a spreading speed $c^*(e)$ in the sense  that, a new mutant or population which is initially confined to a bounded set spreads asymptotically at speed $c^*(e)$ in direction $e$. This falls under the general  Weinberger type given in equation (\ref{op}).

 Note that for the standard Fisher's equation with the drift in one dimension that is
\begin{equation}\label{kppeq}
u_t = d u_{xx}-a u_x+f(u) \ \ \text{for} \ \ (x,t)\in \Omega\times\mathbb R^+,
\end{equation}
the critical domain size for the persistence versus extinction is 
\begin{equation}\label{Lstar} 
L^*:= \frac{2\pi d}{\sqrt{4d f'(0)- a^2}}  ,
\end{equation}
when $\Omega=(0,L)$ and the speeds of propagation to the right and left are
\begin{equation}\label{cstara} 
 c_{\pm}^*(a)=2\sqrt{d f'(0)} \pm a,
 \end{equation}
when $\Omega=\mathbb R$.    For more information regarding the minimal domain size, we refer interested readers to  Lewis et al. \cite{lhl}, Murray and  Sperb \cite{ms},  Pachepsky et al. \cite{plnl},  Speirs and Gurney \cite{sg} and references therein.    The remarkable point is that $c_{\pm}^*(a)$ is a linear function of $a$ and  $L^*$ as a function of $a$ blows up to infinity exactly at roots of $c_{\pm}^*(a)$. For more information on analysis of reaction-diffusion models and on strong connections between persistence criteria and propagation speeds, we refer interested readers to  \cite{lhl, plnl,m1,m2,f}  and references therein.

 \begin{notation}
 Throughout this paper the matrix $I=(\delta_{i,j})_{i,j=1}^n$ stands for the identity matrix,  the matrix $A$ is defined as $A=(a_{i,j})_{i,j=1}^n$ and the vector $ a$ is $ a=(a_i)_{n}$. The matrix $A$ and the vector $ a$ have constant components unless otherwise is stated.   We shall refer to vector fields $a$ with $\div a=0$ as divergence free vector fields. 
The notation  $j_{m,1}$ stands for the first positive zero of the Bessel function $J_m$ for any $m\ge 0$ and $\Gamma$ refers to the Gamma function. 
 \end{notation}

 The organization of the paper is as follows. We investigate how the geometry and size of  the domain $\Omega$ affects persistence vs extinction of the species. We  consider various  types of domain,  including  a $n$-hyperrectangle, a ball of radius $R$, and  a general domain with smooth boundary,  to construct critical domain sizes and extreme volume size (Section \ref{domain}). 
  We then provide applications of the main results to  models for marine reserve, terrestrial reserve, insect pest outbreaks and populations subject to climate change (Section \ref{app}).    At the end,  we provide proofs for our main results and discussions.   

 \section{Geometry of the domain for persistence vs extinction}\label{domain}
 A habitat boundary not only can be considered as a natural consequence of physical features such as rivers, roads, or (for aquatic
systems) shorelines but also it can come from interfaces between different types of ecological
habitats such as forests and grasslands. Boundaries can induce various  effects in population dynamics. They  can affect movement patterns, act as a source of mortality or resource
subsidy, or can function as a unique environment with its own set of rules for population interactions. See Fagan et al. \cite{fcc} for further discussion.   Note that a boundary  can have different effects on different species. For example,  a road may act as a barrier for some species but  as a source of mortality for others. Since a boundary can have different effects on different species, the presence of a boundary can
influence community structure in ways that are not completely obvious from the ways in
which they affect each species, see Cantrell and Cosner \cite{cc,cc01} for more information. 

In this section, we consider the following impulsive reaction-diffusion model on domains with hostile  boundaries to explore persistence versus extinction
 \begin{eqnarray}\label{mainbd}
 \left\{ \begin{array}{lcl}
    u_t = \div(A\nabla u-au)   + f(u)  \quad  \text{for} \ \  (x,t)\in\Omega\times (0,1]  ,    \\
    u(x,t)=0   \quad  \text{for} \ \  (x,t)\in\partial\Omega\times (0,1]    ,  \\
 u(x,0)=g(N_m(x))  \quad  \text{for} \ \ x\in \Omega   , \\
  N_{m+1}(x):=u(x,1)   \quad  \text{for} \ \ x\in \Omega ,
\end{array}\right.
  \end{eqnarray}
where $\Omega$ is a bounded domain in $\mathbb R^n$, $ a$ is a vector in $\mathbb R^n$ and $A$ is a constant positive definite symmetric matrix.

 In what follows we provide critical domain sizes for various domains depending on the geometry of the domain. We shall postpone the proofs of these theorems to the end of this section.  We start with the critical domain size for a $n$-hyperrectangle ($n$-orthotope).

\begin{thm} \label{minDL}
Assume that  $\Omega=[0,L_1]\times\cdots\times[0,L_n]$ where $L_1,\cdots,L_n $ are positive constants, $A=d(\delta_{i,j})_{i.j=1}^m$ and the advection $a$  is a constant vector field in $ \mathbb R^n$. Suppose also that $f$ satisfies (F0)-(F1),  $g$ satisfies (G0)-(G2) and (\ref{conditionfg}) holds.    Then, critical domain dimensions $L_1^*,\cdots, L^*_n$  satisfy
\begin{equation}
\sum_{i=1}^{n} \frac{1}{[L^*_i]^{2}} = \frac{1}{d \pi^2} \left[\ln (e^{f'(0)} g'(0))- \frac{| a|^2}{4d}\right], 
\end{equation}
where $\ln (e^{f'(0)} g'(0))> \frac{| a|^2}{4d}$.  More precisely, when
\begin{equation}
\sum_{i=1}^{n} \frac{1}{[L_i]^{2}} > \frac{1}{d \pi^2} \left[\ln (e^{f'(0)} g'(0))- \frac{| a|^2}{4d}\right],\end{equation}
 and  $g$ satisfies (G0)-(G1) then $N_m(x)$ decays to zero that is  $\lim_{m\to\infty} N_m(x)=0$ and when
 \begin{equation}
\sum_{i=1}^{n} \frac{1}{[L_i]^{2}} < \frac{1}{d \pi^2} \left[\ln (e^{f'(0)} g'(0))- \frac{| a|^2}{4d}\right],\end{equation}
and $g$ satisfies (G0)-(G2)  then $\liminf_{m\to \infty} N_m(x)\ge \bar N(x)$ where  $\bar N(x)$ is a positive equilibrium.   In addition,  critical domain dimensions can be arbitrarily large when $\ln (e^{f'(0)} g'(0))< \frac{| a|^2}{4d}$.

\end{thm}

Suppose that the domain $\Omega$ is a $n$-hypercube. Then the critical domain size is explicitly given by the following corollary.

\begin{cor}\label{hypercubed}
Suppose that assumptions of theorem \ref{minDL} hold. In addition,  let  $L_1=\cdots=L_n=L>0$.   Then the critical domain dimension is
\begin{equation}\label{Lstarn}
L^*:=\left\{ \begin{array}{lcl}
   2\pi d \sqrt \frac{n}{ 4d[f'(0)+\ln(g'(0))] - | a|^2},  &  \ \  \text{if} \ \  4d[f'(0)+\ln(g'(0))] - | a|^2>0,  \\
   \infty,  & \ \   \text{if} \ \  4d[f'(0)+\ln(g'(0))] - | a|^2<0.
\end{array}\right.
\end{equation}
\end{cor}
We now put the above results in a figure to clarify the relationship between the critical domain dimension and the advection in two dimensions,  $n=2$.

\begin{figure}[H]
\begin{center}
\includegraphics[height=1.8in]{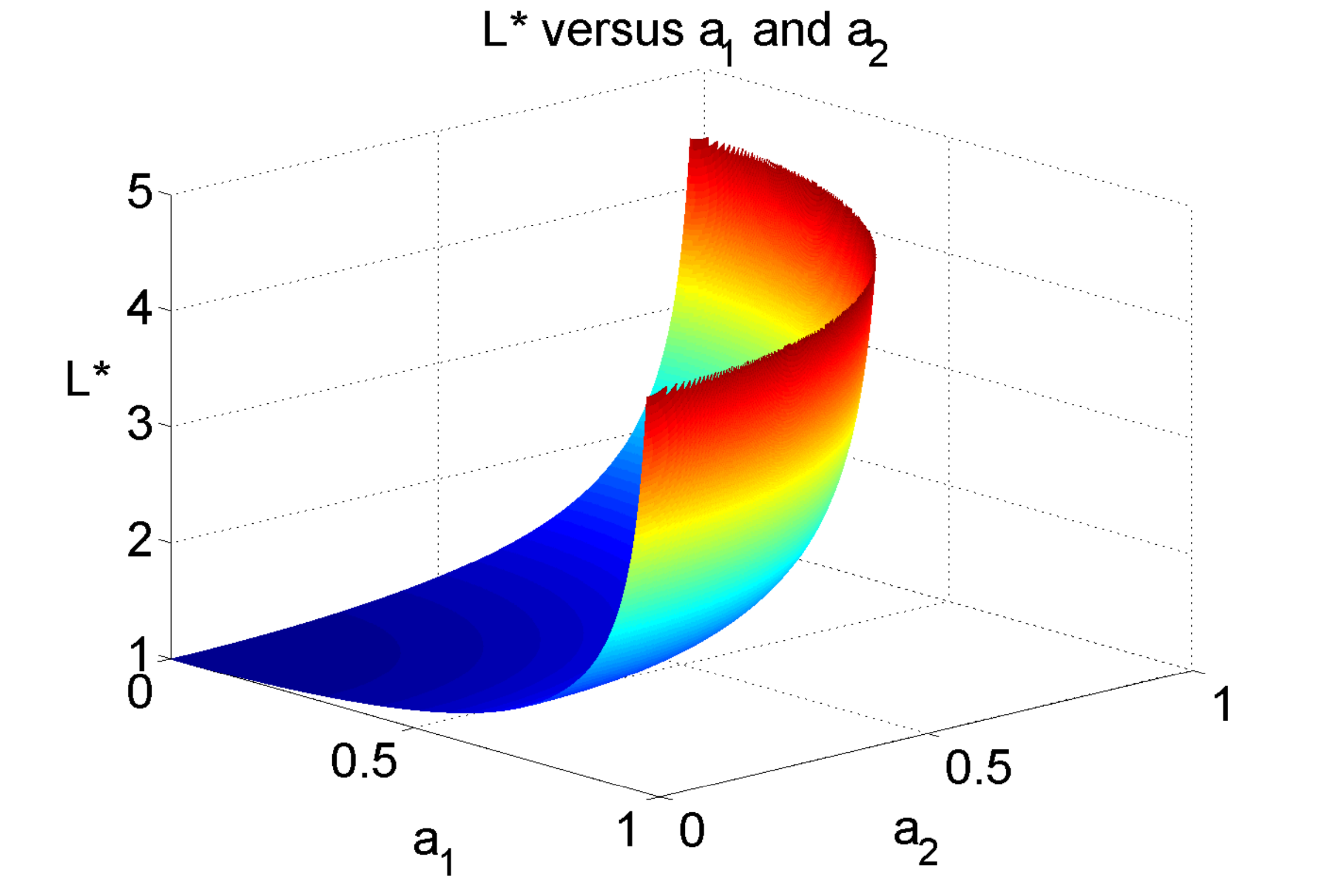}
\vspace*{-10pt}
\caption{Critical domain dimension $L^*$ as a function of $a=(a_1,a_2)$.}
\label{fig:Figure1}
\end{center}
\end{figure}

Note that  for the case $n=1$ and $ a=0$, the minimal domain size $L^*$ was given by Lewis-Li in \cite{LL}.  In addition, one can compare the critical domain dimension (\ref{Lstarn}) for the impulsive reaction-diffusion model (\ref{main}) in higher dimensions to the one given by (\ref{Lstar}) for the classical Fisher's equation in one dimension.   We now provide a similar result for a critical domain radius when the domain is a ball of radius $R$. In Theorem \ref{minDR} and Theorem \ref{minVGeneral},  we let the drift $a$ be a divergence free vector field that is $\div a=0$. 
\begin{thm}\label{minDR}
Suppose that $\Omega=B_R$ where $B_R$ is the ball of radius $R$ and centred at zero, $A=d(\delta_{i,j})_{i,j=1}^m$ and the drift $a$ is a divergence free vector field that is not necessarily constant vector field.   Suppose also that $f$ satisfies (F0)-(F1),  $g$ satisfies (G0)-(G2) and (\ref{conditionfg}) holds.    Then the critical domain radius  is
\begin{equation}
R^*:=j_{n/2-1,1}   \sqrt{\frac{d}{\ln(g'(0)) +f'(0)}}, \end{equation}
in the sense that if $R<R^*$ and $g$ satisfies (G0)-(G1) then $\lim_{m\to\infty} N_m(x)=0$ and if  $R>R^*$ and if $g$ satisfies (G0)-(G2) then  $\liminf_{m\to \infty} N_m(x)\ge \bar N(x)$ where $\bar N(x)$ is a positive equilibrium.
\end{thm}
Theorem \ref{minDR} is a generalized version of the standard island case (two-dimensional space) for discussing the persistence and extinction of a species.

From now on we consider a more general domain $\Omega\subset \mathbb  R^n$ with a smooth boundary.  Our main goal is to provide sufficient conditions on the volume and on the geometry of the domain $\Omega$ for the extinction of populations.     To be able to work with general domains we borrow Schwarz symmetrization rearrangement arguments from mathematical analysis.    In what follows shortly,  we introduce this symmetrization argument and then we apply it to eigenvalue problems associated to (\ref{main}).     We refer interested readers to the lecture notes of Burchard \cite{bur} and references therein for more information.  

Let $\Omega$ be a measurable set of finite volume in $\mathbb R^n$. Its symmetric rearrangement $\Omega^*$ is the following ball centred at zero  whose volume agrees with $\Omega$,
\begin{equation}
 \Omega^*:=\{  x\in\mathbb R^n; \ \   \alpha_n |x|^n < |\Omega| \}.\end{equation}
For any function $\phi \in L^1(\Omega)$, define the distribution function of $\phi$ as
\begin{equation}
 \mu_\phi(t):=|\{x\in\Omega; \phi(x)>t\}|    ,    \end{equation}
for all  $t\in\mathbb R$. The function $\mu_\phi$ is right-continuous and non increasing and as $t\to\infty$  we have $\mu_\phi\to0$. Similarly, as $t\to-\infty$  we have $\mu_\phi\to|\Omega|$. Now for any $x\in\Omega^*\setminus\{0\}$ define
\begin{equation}
 \phi^*(x):=\sup\{t\in\mathbb R, \mu_\phi(t) \ge \alpha_n |x|^n\}.\end{equation}
The function $\phi^*$ is clearly radially symmetric and non increasing in the variable $|x|$. By construction, $\phi^*$ is equimeasurable with $\phi$. In other words,  corresponding level sets of the two functions have the same volume,
\begin{equation}
\mu_\phi(t) = \mu_{\phi^*}(t)   ,   \end{equation}
for all $t\in\mathbb R$.    An essential property of the Schwarz symmetrization is the
following:  if $\phi\in H_0^1(\Omega)$, then $|\phi|^*\in H_0^1(\Omega^*)$ and
\begin{equation}
 ||  |\phi|^* ||_{L^2(\Omega^*)}= || \phi||_{L^2(\Omega)}  , \end{equation}
and
\begin{equation}
|| \nabla |\phi|^* ||_{L^2(\Omega^*)}= || \nabla\phi||_{L^2(\Omega)}.\end{equation}
 One of the main applications of this rearrangement technique is the resolution of optimization problems for the eigenvalues of some second-order elliptic operators on $\Omega$. Let $\lambda_1(\Omega)$ denote the first eigenvalue of the Laplace operator $-\Delta$ with Dirichlet boundary conditions in an open bounded smooth set $\Omega \subset \mathbb R^n$ that is given by the Rayleigh-Ritz formula
  \begin{equation}
 \lambda_1(\Omega) = \inf_{|| \phi||_{L^2(\Omega)} =1} || \nabla\phi||^2_{L^2(\Omega)}.
 \end{equation} 
 It is well-known that $\lambda_1(\Omega) \ge \lambda_1(\Omega^*)$ and that the inequality is strict unless $\Omega$ is a ball.  We refer interested readers to \cite{k,f} and references therein. Since $\lambda_1(\Omega^*)$ can be explicitly computed, this result provides the classical Rayleigh-Faber-Krahn inequality, which states that
\begin{equation}\label{fff}
\lambda_1(\Omega)  \ge \lambda_1(\Omega^*)=\left( \frac{|B_1|}{|\Omega|}\right)^{2/n} j^2_{n/2-1,1} ,\end{equation}
where $j_{m,1}$ is the first positive zero of the Bessel function $J_m$.     More recently, Hamel, Nadirishvili and Russ in \cite{hnr} provided an extension of the above result to the operator 
\begin{equation}
-\div(A\nabla ) + a\cdot \nabla +V, 
\end{equation}
 with Dirichlet boundary conditions where the symmetric matrix field $A$ is in $W^{1,1}(\Omega)$, the vector field (drift) $a:L^1(\Omega)\to\mathbb R^n$ and potential $V$ that is a continuous function in $\bar \Omega$.  Throughout this paper, we call  a matrix  $A$  {\it uniformly elliptic} on $\bar\Omega$ whenever there exists a positive constant $d$ such that for all $x\in\bar \Omega$ and for all $\zeta\in\mathbb R^n$,
\begin{equation}\label{Aelliptic}
A(x) \zeta\cdot \zeta \ge d |\zeta|^2.
\end{equation}

  Consider the following eigenvalue problem,
\begin{eqnarray}\label{eigenA}
 \left\{ \begin{array}{lcl}
   \hfill  - \div(A \nabla \phi) + a \cdot \nabla \phi  &=& \lambda(A, a,\Omega) \phi  \quad  \text{in} \ \ \Omega, \\
   \hfill u & =& 0   \quad  \text{on} \ \  \partial\Omega,
    \end{array}\right.
  \end{eqnarray}
 when $\Omega\in C^{2,\alpha}$ for some $0<\alpha<1$,  the matrix $A$ is uniformly elliptic and $ a :L^\infty(\Omega)\to \mathbb R^n$ where  $|| a||_\infty\le \tau$ for $\tau\ge 0$.  It is shown in \cite{hnr}  that  the first eigenvalue of (\ref{eigenA}) admits the following lower bound 
\begin{equation}\label{lambda1A}
\lambda_1(A, a, \Omega) \ge \lambda_1(d (\delta_{i,j}),\tau e_r,\Omega^*).
\end{equation}
Here $\Omega^*=B_{R}$ is the ball of radius
\begin{equation}
R=\left(\frac{|\Omega|}{\alpha_n}\right)^{1/n},
\end{equation}
for $\alpha_n:=|B_1|=\frac{\pi^{\frac{n}{2}}}{\Gamma(1+\frac{n}{2})}$.  In addition, equality in  (\ref{lambda1A}) holds only when, up to translation, $\Omega=\Omega^*$ and $ a=\tau e_r$ for $e_r=\frac{x}{|x|}$.

 We are now ready to develop sufficient conditions for the extinction of population on various geometric domains.  The remarkable point here is that the volume of the domain $\Omega$ is the key point for the extinction rather than the other aspects of the domain.  The following theorems imply that for some general domain $\Omega$  there is an extreme volume size,  $V_{ex}$, such that  when $|\Omega|<V_{ex}$  extinction must occur for the population living in the habitat.  The following theorem provides an explicit formula for such an extreme volume size for any open bounded domain with a smooth boundary.

\begin{thm}\label{minVGeneral}
Let $A$ be uniformly elliptic,  $f$ satisfy (F0)-(F1),  $g$ satisfy (G0)-(G2) and (\ref{conditionfg}) hold.     Suppose that the vector field $a$ is divergence free and $\Omega\subset \mathbb R^n$ is an open bounded domain with smooth boundary. Then $\lim_{m\to\infty} N_m(x)=0$ for any $x\in\Omega$  when $ |\Omega| < V_{ex}$ for
\begin{equation}\label{Vex1}
V_{ex}:=\frac{1}{\Gamma(1+\frac{n}{2})} \left({\frac{d \pi j^{2}_{\frac{n}{2}-1,1}}{ f'(0)+\ln(g'(0))  } }\right)^{\frac{n}{2}} ,
\end{equation}
and $j_{k,1}$ stands for the first positive root of the Bessel function $J_k$.    

\end{thm}

Theorem \ref{minVGeneral} states that the addition of incompressible flow to the impulse reaction-diffusion problem can only make the critical domain size larger but never smaller. Note that the above theorem, unlike Theorem \ref{minDR} and Theorem \ref{minDL}, only deals with the extinction and not persistence of species. This is due to the fact that, in our proofs, we apply the well-known Rayleigh-Faber-Krahn inequality, see (\ref{fff}), and rearrangement arguments to estimate the first eigenvalue of  elliptic operators. Unfortunately, these estimates are sharp only for the radial domains that is when $\Omega=B_R$. Note also that for the $n$-hyperrectangle,  eigenvalues of elliptic operators are explicitly known. In short, due to a lack of mathematical techniques,  we provide a criteria for the extinction of species and proving persistence under the assumption $ |\Omega| \ge V_{ex}$ seems to be more challenging and remains open.

\begin{cor}\label{coeVex2} In two dimensions,  it is known that $j_{\frac{n}{2}-1,1}=j_{0,1} \approx 2.408$.  Therefore the extreme  volume size is
\begin{equation}\label{Vexn2}
V_{ex}=\frac{d \pi j^2_{0,1}}{f'(0)+\ln(g'(0))}.\end{equation}
\end{cor}

\begin{cor} \label{coeVex3} In three dimensions,  $j_{\frac{n}{2}-1,1}=j_{1/2,1}=\pi$ and $\Gamma(1+\frac{n}{2})=\Gamma(\frac{5}{2})=\frac{3\sqrt \pi}{4}$. Therefore,
\begin{equation}\label{Vexn3}
V_{ex}=\frac{4\pi^4d}{3} \sqrt{  \frac{d }{\left[f'(0)+\ln(g'(0)) \right]^3 } }.
\end{equation}
\end{cor}

The fact that a hyperrectangle (or  $n$-orthotope)   does not have a smooth boundary implies  that Theorem \ref{minVGeneral} is not applied to these types of geometric shapes.  In what follows,  we provide an extreme volume size result for a hyperrectangle  in the light of Theorem \ref{minDL}.

\begin{thm}\label{minVL}
Suppose that $\Omega=[0,L_1]\times\cdots\times[0,L_n]$ and  $A=d (\delta_{i,j})_{i,j=1}^n$ where $L_1,\cdots,L_n $ and $d$ are positive constants and $ a$ is a constant vector field.  Assume that $f$ satisfies (F0)-(F1) and $g$ satisfies (G0)-(G2) and (\ref{conditionfg}) holds.   If $ |\Omega| \le V_{ex}$
for
\begin{equation}
V_{ex}:=\left\{ \begin{array}{lcl}
  \left(   \frac{4 d^2 \pi^2 n}{ 4d[f'(0)+\ln(g'(0))] - | a|^2 } \right)^{n/2}  & \quad  \text{if}  \ \ \ 4d[f'(0)+\ln(g'(0))] >| a|^2,  \\
   \infty,  & \quad  \text{if} \ \  \ 4d[f'(0)+\ln(g'(0))] < | a|^2.
\end{array}\right.
\end{equation}
Then $\lim_{m\to\infty} N_m(x)=0$ for any $x\in\Omega$.
\end{thm}

We shall end this section by briefly explaining the divergence free concept that appeared in  Theorem \ref{minDR} and Theorem \ref{minVGeneral}. The assumption $\div a=0$ implies that the flow field is incompressible.  This would be the case, for example, for water, but not for air.  It also means that the following integral vanishes
\begin{equation}
\int_{\partial\Gamma} a\cdot n \ dS = 0,
\end{equation}
by the divergence theorem.  Here $\Gamma$ is a smooth subset of the domain $\Omega$ (possibly non-proper), and $n$ is the outwardly oriented unit normal vector. The integral quantity which is interpreted as the total flow out of $\Gamma$. This implies that the addition of incompressible flow to the problem cannot make the critical volume size any smaller, although it could make it larger.   Indeed, even in the simple case of one dimension an incompressible flow term of the form $a\equiv C$, where $C$ is a constant, can drive the critical domain size to infinity.

\section{Applications}\label{app}

In this section, we provide applications of main theorems given in the past sections.

\subsection{Marine reserve}\label{mr}
A marine reserve is a marine protected area against fishing and harvesting. Marine reserves could increase species diversity, biomass, and fishery production, etc. within reserve areas,  Lockwood et al. (2002) in \cite{lock}. We use our model here to show the dependence of the critical domain size on the advection flow speed and the mortality rate. Consider the following model
\begin{eqnarray}\label{mainmr}
 \left\{ \begin{array}{lcl}
    u_t = \div(d\nabla u-au)   -\gamma u  \quad  \text{for} \ \  (x,t)\in\Omega\times (0,1], \\
    u(x,t)=0   \quad  \text{for} \ \  (x,t)\in\partial\Omega\times (0,1] , \\
 u(x,0)=g(N_m(x))  \quad  \text{for} \ \ x\in \Omega, \\
  N_{m+1}(x):=u(x,1)   \quad  \text{for} \ \ x\in \Omega ,
\end{array}\right.
  \end{eqnarray}
  where $\gamma$ is the natural mortality rate, and $g$ is either the Beverton-Holt function $g(N_m) =\frac{(1+\lambda)N_m}{1+\lambda N_m}$ or the Ricker function $g(N_m)= N_m e^{r(1-N_m)}$ for $r,\lambda>0$.  The assumption (\ref{conditionfg}) is equivalent to  $\ln(1+\lambda)>\gamma$ or $r>\gamma $ for the Beverton-Holt or Ricker functions, respectively.  Note that the positive equilibrium $N^*$ exists  for  the   Beverton-Holt and Ricker functions and they are of the form 
  \begin{equation}\label{Nstar}
  N^* :=\frac{1}{\lambda} \left[ e^{-\gamma}(1+\lambda)-1 \right] \ \ \ \text{and} \ \ \  N^*: = \frac{r-\gamma}{r},
  \end{equation}
   respectively.  Note also that  the Ricker function is nondecreasing in $[0,N^*] $ when $r\le \gamma+1$ and  the Beverton-Holt function is increasing on the entire $\mathbb R^+$.       We now provide the critical dimension size regarding extinction versus persistence for various geometric shapes.
  \\
  \\
  {\bf Case 1}. Consider the domain $\Omega=[0,L]\times[0,L]$ in two dimensions.    Then Theorem \ref{minDL} implies that the critical dimension for the Beverton-Holt function is
  \begin{equation}\label{c1lb}
  L^*(a)= \frac{\pi \sqrt{2d}}{\sqrt{ \ln(1+\lambda)-\gamma-\frac{|a|^2}{4d}  }} ,
  \end{equation}
  when $ \ln(1+\lambda)>\gamma + \frac{|a|^2}{4d} $  and $L^*(a)$ can be arbitrarily large when  $\gamma<\ln(1+\lambda)<\gamma + \frac{|a|^2}{4d} $.   Similarly, for the Ricker function the critical dimension is
   \begin{equation}\label{c1lr}
  L^*(a)= \frac{\pi \sqrt{2d}}{\sqrt{ r-\gamma-\frac{|a|^2}{4d}  }} ,
  \end{equation}
  when $ \gamma+1\ge r>\gamma + \frac{|a|^2}{4d} $  and $L^*(a)$ can be arbitrarily large when  $r-\gamma $ is close to $\frac{|a|^2}{4d} $.  In other words, when $L<L^*(a)$  the functional sequence $N_m(x)$ satisfies $\lim_{m\to \infty} N_m(x)=0$ that refers to extinction and when $L>L^*(a)$ we have the persistence that is $\lim_{m\to \infty} N_m(x)=\bar N(x)$ for the positive equilibrium $\bar N(x)$.
 \\
  \\
  {\bf Case 2}. Consider the domain $\Omega=B_R$ where $B_R$ is a disk in two dimensions with radius $R$  and $a$ is a divergence free vector field.  Theorem \ref{minDR} implies that the critical value dimension for the Beverton-Holt function is
  \begin{equation}\label{c2rb}
  R^*(\gamma)= j_{0,1} \sqrt{  \frac{d}{\ln(1+\lambda)- \gamma}} ,
  \end{equation}
  and for the Ricker function the critical dimension is
  \begin{equation}\label{c2rr}
  R^*(\gamma)= j_{0,1} \sqrt{  \frac{d}{r- \gamma}} ,
  \end{equation}
  where $j_{0,1}\approx 2.408$ is the first positive root of the Bessel function $J_0$.   In other words, for subcritical radius that is when $R<R^*(\gamma)$ we have  $\lim_{m\to \infty} N_m(x)=0$ and for supercritical radius  that is when  $R>R^*(\gamma)$ we have $\lim_{m\to \infty} N_m(x)=\bar N(x)$ for some positive equilibrium $\bar N(x)$.
    \\
  \\
  {\bf Case 3}. Consider a general domain $\Omega$ in two dimensions with smooth boundary and area $|\Omega|$.  In addition, let  $a$ be a divergence free vector field that is $\div a=0$.  We now apply Theorem \ref{minVGeneral} to find the following  extreme mortality parameter
  \begin{equation}\label{c3rb}
 \gamma_{ex}:=\frac{d j^2_{0,1}}{|\Omega|} -\ln(1+\lambda),
  \end{equation}
  which implies that extinction occurs  for $\gamma< \gamma_{ex}$ and the Beverton-Holt function. Similarly for the Ricker function the extreme mortality parameter $h$ is given by
 \begin{equation}\label{c3rr}
  \gamma_{ex}:=\frac{d j^2_{0,1}}{|\Omega|} -r.
  \end{equation}
  
Note that our main results regarding the critical domain size and the extreme volume size (Section \ref{domain}) are valid in three dimensions as well. 
 \\
  \\
  {\bf Case 4}.  Consider the domain $\Omega=[0,L]\times[0,L]\times[0,L]$ in three dimensions. The corresponding formulae to (\ref{c1lb}) and (\ref{c1lr}) are 
\begin{equation}\label{}
  L^*(a)= \frac{\pi \sqrt{3d}}{\sqrt{ \ln(1+\lambda)-\gamma-\frac{|a|^2}{4d}  }}  \ \ \ \text{and } \ \ \ 
  L^*(a)= \frac{\pi \sqrt{3d}}{\sqrt{ r-\gamma-\frac{|a|^2}{4d}  }} ,
  \end{equation}
for the Beverton-Holt  and  Ricker functions, respectively, under the same assumptions on the parameters.  We now suppose that $\Omega= B_R$ when $B_R$ is a disk in three dimensions with radius $R$  and $a$ is a divergence free vector field. The corresponding formulae to (\ref{c2rb}) and (\ref{c2rr}) are 
\begin{equation}\label{}
  R^*(\gamma)= \pi \sqrt{  \frac{d}{\ln(1+\lambda)- \gamma}} \ \ \ \text{and} \ \ \ 
  R^*(\gamma)= \pi \sqrt{  \frac{d}{r- \gamma}} ,
  \end{equation}
for the Beverton-Holt  and  Ricker functions, respectively.  Lastly,  we consider a general domain $\Omega$ in three dimensions with smooth boundary and area $|\Omega|$. The corresponding extreme mortality parameters to (\ref{c3rb}) and (\ref{c3rr})  are
  \begin{equation}
 \gamma_{ex}:=d \left( \frac{ 4 \pi^4 }{3|\Omega|}\right)^{\frac 2 3} -\ln(1+\lambda) \ \ \ \text{and} \ \ \ 
  \gamma_{ex}:=d \left( \frac{ 4 \pi^4 }{3|\Omega|}\right)^{\frac 2 3} -r, 
  \end{equation}
for the Beverton-Holt  and  Ricker functions, respectively. 

As an application of these results, one can consider the water volume of a marine reserve region for fish. Reserves are being applied commonly to conserve  fish and other populations under threat, see \cite{bae,ms}. We know that, because water is incompressible, the currents will give a divergence free vector field in three dimensions.  Then Theorem  \ref{minVGeneral} implies that ocean currents could not improve the persistence of the species, under the assumption of hostile boundary conditions at the exterior of the reserve.

\subsection{Terrestrial reserve}
A terrestrial reserve is a terrestrial protected area for conservation and economic purposes. Terrestrial reserves  conserve biodiversity and protect threatened and endangered species from hunting \cite{fah}. We use our model here to show the dependence of terrestrial species persistence on the diffusion rate.
Consider
\begin{eqnarray}\label{mainmr}
 \left\{ \begin{array}{lcl}
    u_t = d \Delta u  -\gamma u  \quad  \text{for} \ \  (x,t)\in\Omega\times (0,1],  \\
    u(x,t)=0   \quad  \text{for} \ \  (x,t)\in\partial\Omega\times (0,1] , \\
 u(x,0)=g(N_m(x))  \quad  \text{for} \ \ x\in \Omega  , \\
  N_{m+1}(x):=u(x,1)   \quad  \text{for} \ \ x\in \Omega ,
\end{array}\right.
  \end{eqnarray}
  where $g$ is the Beverton-Holt function $g(N_m) =\frac{(1+\lambda)N_m}{1+\lambda N_m}$.  Just like the previous Section \ref{mr} the assumption (\ref{conditionfg}) is equivalent to  $\ln(1+\lambda)>\gamma$.
  \\
  \\
  {\bf Case 1}. Consider the domain $\Omega=[0,L_1]\times[0,L_2]$ in two dimensions where $L_1,L_2$ are positive.     Then a direct consequence of  Theorem \ref{minDL} is that when the diffusion is greater than the following critical diffusion
  \begin{equation}\label{dextL}
  d^*:= \frac{1}{\pi^2} [\ln(1+\lambda)-\gamma]\frac{L_1^2 L^2_2}{L_1^2+L_2^2},
  \end{equation}
  that is when $d>d^*$ we have $\lim_{m\to \infty} N_m(x)=0$ that refers to extinction.  In addition, when $d<d^*$ we have  $\lim_{m\to \infty} N_m(x)=\bar N(x)$ for some positive equilibrium $\bar N(x)$.
   \\
  \\
  {\bf Case 2}. We now consider a bounded domain $\Omega$ in two dimensions with a smooth boundary and area $\Omega$.  Theorem \ref{minVGeneral} implies that  when the diffusion is greater than the following extreme diffusion
  \begin{equation}\label{dext}
  d_{ex}:=  \frac{|\Omega|}{\pi j_{0,1}^2}    [\ln(1+\lambda)-\gamma],
  \end{equation}
that is when $d>d_{ex}$ then population must go extinct.   Note that for the case of $\Omega=B_R$, Theorem \ref{minDR} implies that  the extreme diffusion given by (\ref{dext}) is the critical diffusion
  \begin{equation}
  d^*:=  \frac{R^2}{j_{0,1}^2}    [\ln(1+\lambda)-\gamma] ,
  \end{equation}
  meaning that when $d>d^*$ then $\lim_{m\to \infty} N_m(x)=0$ and for $d<d^*$ we have  $\lim_{m\to \infty} N_m(x)=\bar N(x)$ for some positive equilibrium $\bar N(x)$.   In Case 1,  we considered a rectangular domain that does not have smooth boundaries unlike the domain in Case 2.  Comparing $d^*$ and $d_{ex}$ given by  (\ref{dextL}) and (\ref{dext}) one sees that $d^*<d_{ex}$  when $\Omega=[0,L_1]\times[0,L_2]$ and $L_1=L_2$, since  $\pi j_{0,1}^2 \approx 18.21$ and $ 2\pi^2=19.73$.

\subsection{Insect pest outbreaks}
Insect pest outbreaks are a historic problem in agriculture and can have long-lasting effects. It is may be  necessary to control insect pests in order to maximize crop production \cite{tc}. We apply our model here to study the dependence of insect pest extirpation on the removal rate. Consider
\begin{eqnarray}\label{mainmr}
 \left\{ \begin{array}{lcl}
    u_t = d \Delta u + r(1 -u)u  \quad  \text{for} \ \  (x,t)\in\Omega\times (0,1], \\
    u(x,t)=0   \quad  \text{for} \ \  (x,t)\in\partial\Omega\times (0,1] ,\\
 u(x,0)=g(N_m(x))  \quad  \text{for} \ \ x\in \Omega  , \\
  N_{m+1}(x):=u(x,1)   \quad  \text{for} \ \ x\in \Omega,
\end{array}\right.
  \end{eqnarray}
  where  $g(N_m) =(1-s) N_m$ and $0<s<1$ is the surviving fraction of pests that contribute to the population a year later.   The assumption (\ref{conditionfg}) is equivalent to  
  \begin{equation}\label{ers}
  e^r (1-s) > 1.
  \end{equation} 
  Assuming that (\ref{ers}) holds,  straightforward calculations show that the positive equilibrium $N^*$ that solves (\ref{nstar}) is of the form 
  \begin{equation}\label{Nstarr}
  N^*:= \frac{ (1-s)e^r -1 }{(1-s) (e^r-1)},
   \end{equation} 
  and $g$ is clearly increasing in $[0,N^*]$. 
  \\
  \\
  {\bf Case 1}.  Consider the rectangular domain $\Omega=[0,L_1]\times[0,L_2]$ in two dimensions for positive $L_1,L_2$. We now apply Theorem \ref{minDL} to conclude that  the critical value for the parameter $s$ regarding extinction versus persistence is
   \begin{equation}
s^*:= 1-e^{d\pi^2       \left[ \frac{L_1^2+L_2^2}{L_1^2L_2^2} \right]-r} .
    \end{equation}
 More precisely, for $s>s^*$ we have $\lim_{m\to \infty} N_m(x)=0$ and for $s<s^*$ we have  $\lim_{m\to \infty} N_m(x)=\bar N(x)$ for the positive equilibrium $\bar N(x)$.
   \\
  \\
  {\bf Case 2}. Consider a slightly more general bounded domain $\Omega$ in two dimensions with a smooth boundary and area $|\Omega|$.  Here we introduce an extreme value for the      parameter $s$ as
  \begin{equation}
  s_{ex}:= 1-e^{\frac{d \pi j_{0,1}^2}{|\Omega|}-r} ,
    \end{equation}
meaning that when $s$ is larger than $s_{ex}$ the population must go extinct.

\subsection{Populations in moving habitats arising from climate change}

We end this part with mentioning that  (\ref{mainmr}) can be applied to study population subject to  climate change. Climate change, especially global warming, has greatly changed the distribution and habitats of biological species. Uncovering the potential impact of climate change on biota is an important task for modelers \cite{pl}.   We consider a rectangular domain that is $\Omega=[0,L_1]\times[0,L_2]$ moving in the positive $x$-axis direction at speed $c$. Outside this domain conditions are hostile to population growth, while inside the domain there is random media, mortality and periodic reproduction. Using the approach of \cite{pl} this problem is transformed to a related problem on a stationary domain. Consider (\ref{mainmr}) where $g$ is  the Beverton-Holt function  and the vector field $a$ is non zero in $x$-axis direction that is $a=(-c,0)$.    Theorem \ref{minDL} implies that for
  \begin{equation}\label{L1L2}
  \frac{1}{L_1^2} + \frac{1}{L_2^2}<  \frac{1}{d\pi^2} \left[\ln(1+\lambda)-\gamma - \frac{c^2}{4d}\right],
  \end{equation}
when   $\ln(1+\lambda)>\gamma + \frac{c^2}{4d}$ we have  $\lim_{m\to \infty} N_m(x)=\bar N(x)$ for some positive equilibrium $\bar N(x)$ which refers to persistence of population.  In other words, (\ref{L1L2}) yields the parameter $c$ must be bounded by
\begin{equation}\label{cL}
  c^2 < 4d[\ln(1+\lambda)- \gamma ] - (2d\pi)^2 \left[\frac{L_1^2+L_2^2}{L_1^2L_2^2} \right] <4d[\ln(1+\lambda)- \gamma ] .
\end{equation}
   This implies that  persistence and ability
to propagate should be closely connected.  This is the case from a biological perspective as well. For example, if a population
cannot propagate upstream but is washed downstream,  it will not persist.   We refer interested readers to Speirs and  Gurney \cite{sg} and Pachepsky et al. \cite{plnl} for similar arguments regarding  Fisher's equation with advection.

\section{Discussion} \label{secdis}
We examined impulsive reaction-diffusion equation models for species with distinct reproductive and dispersal stages on domains $\Omega\subset \mathbb R^n$ when $ n\ge1$ with diverse geometrical structures.   Unlike standard partial differential equation models,  study of impulsive reaction-diffusion models requires a simultaneous analysis of the differential equation and the recurrence relation.  This fundamental fact rules out certain standard mathematical analysis theories for analyzing solutions of these type models,  but it  opens  up various ways to apply the model.   These models can be considered as a description for a continuously growing and dispersing population with pulse harvesting and a population with individuals immobile during the winter.  As a domain $\Omega$ for the model, we considered bounded sets in $\mathbb R^n$, including  convex and non convex domains with smooth or non smooth boundaries, and also the entire space  $\mathbb R^n$. Since the geometry of the domain has  tremendous impacts on the solutions of the model, we consider various type domain to study the qualitative properties of the solutions.  We refer interested readers to \cite{fcc,cc01,cc} and references therein regarding {\it how habitat edges change species interactions}. 

On bounded rectangular and circular domains, we provided critical domain sizes regarding  persistence versus extinction of populations in any space dimension $n \ge 1$.   In order  to find the critical domain sizes we used the fact that the first eigenpairs  of the  Laplacian operator for such domains can be computed explicitly. Note that for a general bounded domain in $\mathbb R^n$, with smooth boundaries,  the spectrum of the Laplacian operator is not known explicitly. However, various estimates are known for the eigenpairs.  As a matter of fact, study of the eigenpairs of  differential operators  is one of the oldest problems in the field of  mathematical analysis and partial differential equations, see Faber \cite{f}, Krahn \cite{k}, P\'{o}lya \cite{p} and references therein.  

  We applied several mathematical analysis methods such as Schwarz symmetrization rearrangement arguments,  the classical Rayleigh-Faber-Krahn inequality and  lower bounds for the spectrum of uniformly elliptic operators with Dirichlet boundary conditions given by Li-Yau \cite{ly} to compute a novel quantity called extreme volume size.   Whenever $|\Omega|$ falls below the  extreme volume size  the species must driven extinct,    regardless of the geometry of the domain. In other words, the extreme volume size provides a lower bound (a necessary condition) for the persistence of population. In the context of biological sciences this can provide a partial answer to questions about {\it how much habitat is enough?} (see also Fahrig \cite{fah} and references therein).  We believe that this opens up new directions of research in both mathematics and sciences.

 In this paper we presented applications of our main results in two and three space dimensions to certain biological reaction-diffusion models regarding  marine reserve, terrestrial reserve, insect pest outbreaks and  population subject to climate change.  Note that when $n=1$ our results recovers the results provided by  Lewis and Li \cite{LL} and when $g(N)=N$ and $A=d(\delta_{i,j})_{i,j=1}^n$ our results on both bounded and unbounded domains coincide with the ones for the standard Fisher-KPP equation.    Throughout this paper we assumed that the diffusion matrix and, for the most parts, the advection vector field contain constant components. Study of the influence of non constant advection and non constant diffusion on persistence and  extinction properties can be considered along the work of Hamel and Nadirashvili \cite{hn} and Hamel, Nadirashvili and Russ \cite{hnr}.

\section{Proofs}\label{secproofs}

In this section, we provide mathematical proofs for our main results.  

\subsection{Proofs of Theorem \ref{minDL} and Theorem \ref{minDR}}
We start with the proof of Theorem \ref{minDL}. Proof of Theorem \ref{minDR} is very similar. For both of the proofs, the following  technical lemma plays a fundamental role. We shall omit the proof of the lemma since it is standard in this context. 
 
\begin{lemma}\label{eigenvaluesimple}
Suppose that $\phi_1,\lambda_1$ are the first eigenvalue and the first eigenfunction of (\ref{eigenA}). Then
 \begin{eqnarray}\label{eigenA2}
&&\lambda_1(d I, a,[0,L_1]\times\cdots\times[0,L_n])= \frac{| a|^2}{4d} +d \pi^2 \left( \frac{1}{L^2_1}+ \cdots+\frac{1}{L^2_n} \right), \\&& \lambda_1(d I,0,B_R)=   j^2_{n/2-1,1} R^{-2} d,
  \end{eqnarray}
  where $j_{m,1}$ is the first positive zero of the Bessel function $J_m$.
\end{lemma}

\noindent {\it Proof of Theorem \ref{minDL}}. Suppose that  the recurrence relation $\bar N_m(x)$ is a solution of  the following linearized problem
 \begin{eqnarray}\label{mainbdlin}
 \left\{ \begin{array}{lcl}
    u_t + a \cdot \nabla u = \div(A \nabla u)   + f'(0)u  \quad  \text{for} \ \  (x,t)\in\Omega\times \mathbb R^+, \\
    u(x,t)=0   \quad  \text{for} \ \  (x,t)\in\partial\Omega\times\mathbb R^+, \\
 u(x,0)=g'(0)\bar N_m(x)  \quad  \text{for} \ \ x\in \Omega,\\
  \bar N_{m+1}(x):=u(x,1)   \quad  \text{for} \ \ x\in \Omega.
\end{array}\right.
  \end{eqnarray}
Let 
\begin{equation}\label{uxtt}
u(x,t)=c g'(0) e^{\lambda t} \phi(x),
\end{equation} 
 be a solution for the above linear problem (\ref{mainbdlin}) for some function $\phi$ and constant $c$.  It is straightforward to show that $\phi$ satisfies 
\begin{equation}\label{Aeigen}
-\div(A \nabla \phi)+ a \cdot \nabla \phi=(f'(0)-\lambda)\phi \ \ \text{in } \ \ \Omega. 
\end{equation}\label{limit}
 Note that the initial condition is $u(x,0)=g'(0)\bar N_m(x)=cg'(0) \phi(x) $. This implies that  
 \begin{equation}
  \bar N_{m+1}(x)=u(x,1)=c g'(0)  e^{\lambda} \phi(x) = g'(0)  e^{\lambda} \bar N_m(x).
  \end{equation}
       In the light of (\ref{Aeigen}),  we consider the following Dirichlet eigenvalue problem
 \begin{eqnarray}\label{eigen}
 \left\{ \begin{array}{lcl}
   \hfill -\div(A \nabla \phi)+ a \cdot \nabla \phi&=&\lambda \phi  \quad  \text{in} \ \  \Omega ,\\
 \hfill \phi&=&0   \quad  \text{on} \ \  \partial\Omega ,\\
\hfill  \phi&>&0   \quad  \text{in} \ \  \Omega , \\
\end{array}\right.
  \end{eqnarray}
and let the pair $(\phi_1,\lambda_1(A, a,\Omega))$ be the first eigenpair of this problem. Setting $\lambda:= f'(0)-\lambda_1(A, a,\Omega)$ in (\ref{uxtt}) we get 
\begin{equation}
u(x,t)=c g'(0) e^{(f'(0)-\lambda_1(A, a,\Omega)) t} \phi_1(x).
\end{equation}
This implies that
\begin{equation}\label{}
\bar N_{m+1}(x)= c \left(g'(0) e^{f'(0)-\lambda_1(A, a,\Omega)} \right)^m \phi_1(x).
\end{equation}
Therefore, we get
\begin{equation}\label{limit}
\lim_{m\to\infty} \bar N_m(x)=0,
\end{equation}
when $g'(0) e^{f'(0)-\lambda_1(A, a,\Omega)}<1$,  that is when
\begin{equation}\label{lambdaA}
\lambda_1(A, a,\Omega) > \ln (e^{f'(0)} g'(0)).
\end{equation}
Suppose that $u(x,0)=N_0(x)$ is an initial value for the original nonlinear problem (\ref{main}). One can choose a sufficiently large $c$ such that $N_0(x) \le \bar N_0(x)$.  Applying standard comparison theorems, together with the fact that $g$ is linear, and induction arguments we obtain  $N_m(x)\le \bar N_m(x)$ for all $m\ge 0$.  This implies that 
\begin{equation}\label{limitN}
\lim_{m\to\infty} N_m(x)=0,
\end{equation}
when (\ref{lambdaA}) holds.  We now analyze inequality (\ref{lambdaA}) that involves the first eigenvalue of  (\ref{eigen}).   Applying Lemma \ref{eigenvaluesimple} to (\ref{eigen}) when $A=d (\delta_{i,j})_{i,j=1}^m$, $\Omega=[0,L_1]\times\cdots\times[0,L_n]$ we obtain 
\begin{equation}\label{lambdadI}
\lambda_1(dI, a,\Omega)=\frac{| a|^2}{4d} +d \pi^2 \left(\sum_{i=1}^{n} L_i^{-2}\right)   .
\end{equation}
Equating this and  (\ref{lambdaA}) we end up with 
\begin{equation}\label{sumL}
\sum_{i=1}^{n} L_i^{-2} > \frac{1}{4 d^2 \pi^2} \left[ 4d [\ln (e^{f'(0)} g'(0))]- |a|^2 \right].
\end{equation}
This proves the first part of the theorem.   To show the second part we suppose that $g$ satisfies (G0)-(G2) and we  analyze the eigenvalue problem (\ref{eigen}).   Suppose that $g'(0) e^{f'(0)-\lambda_1(A, a,\Omega)}>1$ that is when the inequality (\ref{lambdaA}) is reversed, namely, 
\begin{equation}
\lambda_1(A, a,\Omega)  < \ln(g'(0)) +f'(0)  . 
\end{equation}
 Let $\un{\lambda}>\lambda_1$ and $\un{g}$ such that $\un g e^{f'(0) - \un \lambda}>1$.  Since $g$ satisfies (G0)-(G2),  there exist a $\tilde M$ such that  $g(N)\ge g'(0)N-h(N)$ for $0<N< \tilde M$. The fact that $h$ is differentiable and $h(0)=h'(0)=0$ implies that there exists a positive constant $\delta$ such that for $N<\delta$ the fraction $\frac{h(N)}{N}$ can be sufficiently small, let's say 
\begin{equation}\label{hnn}
 \frac{h(N)}{N}<\min\{g'(0)-\un g, \un \lambda -\lambda_1\}.
 \end{equation}
    Now define 
 \begin{equation}
 \un u(x,t)= \epsilon \un g e^{(f'(0)-\un \lambda)t} \phi_1(x) .
 \end{equation}  
 Note that $g(\un u(x,t)) \ge g'(0) \un u(x,t) - h(\un u(x,t))$. For sufficiently small $\epsilon$ and $t\in (0,1]$, we get 
   \begin{eqnarray}\label{}
\frac{g(\un u(x,t))}{\un u(x,t)} \ge \un g + (g'(0)-\un g) - \frac{h(\un u(x,t))}{\un u(x,t)} \ge \un g ,
  \end{eqnarray}
where we have used the fact that $\frac{h(\un u(x,t))}{\un u(x,t)} <g'(0)-\un g$ when $t\in (0,1]$.  This implies that 
\begin{equation} \label{gunder}
g(\un u(x,t) ) \ge \un g \  \un u(x,t) . 
\end{equation}
We now show that $\un u$ is a subsolution for the partial differential equation given in (\ref{main}) when $t\in (0,1]$  that is 
 \begin{eqnarray}\label{subp}
 \left\{ \begin{array}{lcl}
   \hfill u_t=\div(A\nabla u)+a \cdot \nabla u+ f(u)   \quad  \text{for} \ \  (x,t)\in\Omega\times \mathbb R^+ ,\\
    u(x,t)=0   \quad  \text{for} \ \  (x,t)\in\partial\Omega\times\mathbb R^+.
\end{array}\right.
  \end{eqnarray}
Since $f$ satisfies (F0)-(F1), we get 
 \begin{equation}
f(\un u) \ge  f'(0) \un u - h(\un u  ) =   \epsilon \un g e^{(f'(0)-\un \lambda)t}  \phi_1f'(0) - h(\un u  ) .
\end{equation}
 This implies that 
    \begin{eqnarray}\label{subu}
 && \un u_t -\div(A\nabla \un u)-a \cdot \nabla \un u- f(\un u) 
 \\&\le& \epsilon \un g e^{(f'(0)-\un \lambda)t} \left[(f'(0)-\un \lambda)\phi_1 - \div(A\nabla \phi_1) +a\cdot \nabla \phi_1 - f'(0) \phi_1\right]+ h(\un u)
\nonumber \\&=& \epsilon \un g e^{(f'(0)-\un \lambda)t} \left[- \div(A\nabla \phi_1) +a\cdot \nabla \phi_1 - \lambda_1 \phi_1\right] 
\nonumber \\&&\nonumber+ \epsilon \un g e^{(f'(0)-\un \lambda)t}  \left[  \lambda_1-\un\lambda \right] \phi_1 + h(\un u) .
  \end{eqnarray}
Note that $\lambda_1$ is the first eigenvalue of (\ref{eigen}) that yields 
\begin{equation}
 -\div(A\nabla \phi_1) +a\cdot \nabla \phi_1 - \lambda_1 \phi_1=0.
\end{equation}
 Applying this in (\ref{subu}) we get 
 \begin{eqnarray}\label{subu11}
 && \un u_t -\div(A\nabla \un u)-a \cdot \nabla \un u- f(\un u) 
 \\&\le &\nonumber \epsilon \un g e^{(f'(0)-\un \lambda)t}  \left[  \lambda_1-\un\lambda \right] \phi_1 + h(\un u)   
 \\&= &\nonumber \un u \left[   \lambda_1-\un\lambda + \frac{h(\un u)}{\un u}  \right]. 
 \end{eqnarray}
We now apply (\ref{hnn}) to conclude that 
\begin{equation}\label{subsolu}
\un u_t -\div(A\nabla \un u)-a \cdot \nabla \un u- f(\un u) \le 0.
\end{equation} 
This proves our claim. Suppose that $N_m(x)$ is a solution of (\ref{main}) then 
\bg
N_{m+1}(x)= Q[g(N_m)] (x),
\ed
where the operator $Q$ maps $u(x,0)$ to $u(x,1)$. To be more precise, let $u(x,t)$ solve 
\begin{eqnarray}\label{Qubar}
 \ \ \ \ \  \left\{ \begin{array}{lcl}
    u_t  = \div(A\nabla u-a u)   + f(u)  \quad  \text{for} \ \  (x,t)\in\Omega\times (0,1] \\
 u(x,0)=u_0(x)  \quad  \text{for} \ \ x\in \Omega\\
\end{array}\right.
  \end{eqnarray}
then $u(x,1)=Q[u_0](x)$. It is straightforward to see that $Q$ is a monotone operator due to comparison principles. Now, let $M_0(x)=\epsilon \phi_1(x)$ and 
\bg M_{m+1} (x)=Q[g(M_m)](x). 
\ed
 We now apply the fact that $\un u$ is a subsolution of (\ref{Qubar}) together with  (\ref{subsolu}) and (\ref{gunder}) and the monotonicity of the operator $Q$ to conclude  
 \begin{equation}
 M_1(x)=Q[g(M_0)](x) \ge Q[\un g \cdot  M_0](x) \ge \un u(x,1) \ge M_0(x). 
 \end{equation}
 An induction argument implies that 
  \begin{equation}
M_{m+1}(x) \ge M_m (x)  \ \ \text{for} \ \ m\ge 0. 
  \end{equation}
Note that $g(0)=0$ and $Q$ is a continuous operator defined on the range of function $g$.  This implies that for a sufficiently small $\epsilon>0$, we have $M_0(x)=\epsilon \phi_1(x) \le M_1(x)=Q[g(\epsilon \phi_1(x)   )](x) \le N^*$ where $N^*$ is a positive equilibrium for (\ref{equib}). Therefore, 
  \begin{equation}\label{MmN}
M_{m}(x) \le M_{m+1} (x)  \le N^*.  
   \end{equation}
The sequence of function $M_m(x)$ is bounded increasing that must be convergent to $N(x)$. Note also that for sufficiently small $\epsilon$ one can show that  $M_0(x)=\epsilon \phi_1(x) \le Q[g(u_0)] (x)=Q[g(N_0)] (x)=N_1(x)$. We refer interested readers to \cite{jst} for a similar argument.   Finally from comparison principles we have $M_{m+1}(x) \le N_m(x)$.    This and (\ref{MmN}) imply that $\liminf_{m\to \infty} N_m(x) \ge \liminf_{m\to \infty} M_m(x) =N(x)$. This completes the proof.

\hfill $ \Box$

\noindent {\it Proof of Corollary \ref{hypercubed}}.  Assuming that $L_i=L>0$ for any $1\le i\le n$, from (\ref{sumL}) we get the following
\begin{equation}\label{sumLS}
n L^{-2} > \frac{1}{4 d^2 \pi^2} \left( 4d [\ln (g'(0)) +f'(0)]- |a|^2 \right).
\end{equation}
This completes the proof.

\hfill $ \Box$

\noindent {\it Proof of Theorem  \ref{minDR}}.  Suppose that the recurrence relation $\bar N_m(x)$ is a solution of the linearized problem (\ref{mainbdlin}) when $\div a=0$ and $A=d (a_{i,j})_{i,j=1}^m$.    Assume that the following function $u(x,t)$ is a solution for the linear problem (\ref{mainbdlin}), 
\begin{equation}\label{uxttR}
u(x,t)=c g'(0) e^{\lambda t} \phi(x),
\end{equation} 
 for some constant $c$ and where $\phi$ satisfies 
\begin{equation}\label{AeigenR}
-d\Delta \phi+ a \cdot \nabla \phi=(f'(0)-\lambda)\phi \ \ \text{in } \ \ B_R. 
\end{equation}\label{limit}
Suppose that  $\lambda_1$ is the first eigenvalue of the following eigenvalue problem with Dirichlet boundary conditions
 \begin{eqnarray}\label{eigenB}
 \left\{ \begin{array}{lcl}
   \hfill -d\Delta \phi+ a \cdot \nabla \phi&=&\lambda (dI,a,B_R) \phi  \quad  \text{in} \ \  B_R \\
 \hfill \phi&=&0   \quad  \text{on} \ \  \partial B_R \\
\hfill  \phi&>&0   \quad  \text{in} \ \  B_R \\
\end{array}\right.
  \end{eqnarray}
  Under similar arguments as in the proof of Theorem \ref{minDL}, whenever 
  \begin{equation}\label{lambda1R}
\lambda_1(dI, a,B_R) > \ln (e^{f'(0)} g'(0)), 
\end{equation}
  we  conclude the following decay 
 \begin{equation}\label{limitBR}
\lim_{m\to\infty} N_m(x)= \lim_{m\to\infty} \bar N_m(x)= 0  \ \ \text{in} \ \ B_R, 
\end{equation}
and otherwise we get 
 \begin{equation}\label{limitBRR}
\lim_{m\to\infty} N_m(x) \ge \lim_{m\to\infty} \bar N_m(x) \ge N(x) \ \ \text{in} \ \ B_R. 
\end{equation}
Therefore, we only need to discuss the magnitude of $\lambda_1$.  Lemma \ref{eigenvaluesimple} implies that 
  \begin{equation}\label{lambda1BRUp}
 \lambda_1(dI,0,B_R) =d \left( \frac{|B_1|}{|B_R|}\right)^{2/n} j^2_{n/2-1,1},
   \end{equation}
where $j_{n/2-1,1}$ is the first positive zero of the Bessel function $J_{n/2-1}$. This and the fact that $a$ is a divergence free vector field implies that 
\begin{equation}\label{lambda1BR}
\lambda_1(dI,a,B_R) \ge d R^{-2} j^2_{n/2-1,1} .
\end{equation}
 Combining (\ref{lambda1BR}) and (\ref{lambda1R}) provides proofs for (\ref{limitBR}) and (\ref{limitBRR})  when $R<R^*$ and $R>R^*$, respectively,  where
\begin{equation}
R^{*} := \sqrt{\frac{d \  j^2_{n/2-1,1}}{ \ln (g'(0)) +f'(0) }}.
\end{equation}
This completes the proof.

\hfill $ \Box$

\subsection{Proofs of Theorem \ref{minVGeneral} and Theorem \ref{minVL}}

This part is dedicated to proofs of theorems regarding the extreme volume size for both a general domain $\Omega\subset \mathbb R^n$ with a smooth boundary, Theorem \ref{minVGeneral}, and also for a $n$-hyperrectangle domain $\Omega=[0,L_1]\times \cdots \times [0,L_n]$ with non smooth boundaries,  Theorem \ref{minVL}.   
\\
\\
\noindent {\it Proof of Theorem  \ref{minVGeneral}}.   With the same reasoning as in the proof of Theorem \ref{minDL}, the limit tends to zero that is 
\begin{equation}\label{limitomega}
\lim_{m\to\infty} N_m(x)= \lim_{m\to\infty} \bar N_m(x)= 0  \ \ \text{in} \ \ \Omega, 
\end{equation}
 when 
\begin{equation}\label{lambdaVex}
\lambda_1(A, a,\Omega) > \ln (e^{f'(0)} g'(0)). 
\end{equation}
Here $\lambda_1(A, a,\Omega)$ stands for the first eigenvalue of
\begin{eqnarray}\label{eigenVex}
 \left\{ \begin{array}{lcl}
 \hfill -\div(A \nabla \phi)+ a \cdot \nabla \phi&=&\lambda_1(A, a,\Omega) \phi  \quad  \text{in} \ \  \Omega, \\
 \hfill \phi&=&0   \quad  \text{on} \ \  \partial\Omega, \\
\hfill  \phi&>&0   \quad  \text{in} \ \  \Omega. 
\end{array}\right.
  \end{eqnarray}
From the fact that the vector field $ a$ is divergence free (in the sense of distributions),  we have $ \lambda_1(A, a, \Omega) \ge \lambda_1(A,0,\Omega)$.  This can be seen by multiplying  (\ref{eigenVex}) with $\phi$ and integrating by parts over $\Omega$ and using the fact that $\div a=0$.  We now apply the generalized Rayleigh-Faber-Krahn inequality  (\ref{lambda1A})  that is
\begin{equation}\label{lambda1Aa}
 \lambda_1(A, a, \Omega) \ge  \lambda_1(dI,0,\Omega^*)  ,
 \end{equation}
where $\Omega^*$ is the ball $B_{R}\subset \mathbb R^n$ for the radius $R=\left(\frac{|\Omega|}{|B_1|}\right)^{1/n}$ where $|B_1|$ is the volume of the unit ball that is $|B_1|=\frac{\pi^{\frac{n}{2}}}{\Gamma(1+\frac{n}{2})}$. Note that Lemma \ref{eigenvaluesimple} implies that 
 \begin{equation}\label{lambda1*}
 \lambda_1(dI,0,\Omega^*)=j^2_{n/2-1,1} R^{-2} d   .
 \end{equation}
 Combining  (\ref{lambdaVex}), (\ref{lambda1Aa}) and (\ref{lambda1*}) shows that whenever
 \begin{equation}\label{jnR}
j^2_{n/2-1,1} R^{-2} d> f'(0) +\ln(g'(0)),
\end{equation}
the decay estimate (\ref{limitomega}) holds. Note that  the radius $R$ in (\ref{jnR}) is $\left(\frac{|\Omega|}{|B_1|}\right)^{1/n}$.  This implies that for any domain with the volume $|\Omega|$ that is less than
\begin{equation}\label{}
|\Omega|< |B_1| \left(\frac{d\ j^2_{n/2-1,1} }{f'(0) +\ln(g'(0))}\right)^{\frac{n}{2}} ,
\end{equation}
the population must go extinct. This completes the proof.

\hfill $ \Box$

We now provide a proof for Theorem  \ref{minVL} that is in regards to $n$-hyperrectangle domain $$\Omega=[0,L_1]\times \cdots \times [0,L_n]$$ with non smooth boundaries. 
\\
\\
\noindent {\it Proof of Theorem  \ref{minVL}}. The idea of proof is very similar to the ones provided in  proofs of Theorem \ref{minVGeneral} and Theorem \ref{minDL}.      The decay of the $N_m(x)$ to zero as $m$ goes to infinity  refers to  inequality (\ref{sumL}) that is 
\begin{equation}\label{L1L2Ln}
\frac{1}{L^2_1}+ \cdots+\frac{1}{L^2_n} > \frac{1}{d \pi^2} \left[  f'(0) +\ln(g'(0)) -\frac{|a|^2}{4d} \right].
\end{equation}
Note that when the right-hand side of the above inequality is nonpositive then this is valid for any $L_i$ for $1\le i\le n$. So, we assume that $ f'(0) +\ln(g'(0)) -\frac{|a|^2}{4d} >0$.

On the other hand, the following inequality of arithmetic and geometric means hold
\begin{equation}\label{meanL1L2}
\frac{1}{L^2_1}+ \cdots+\frac{1}{L^2_n} \ge n \sqrt[n]{\frac{1}{L^2_1} \cdots\frac{1}{L^2_n}} = n \left(\frac{1}{L_1\cdots L_n}\right)^{\frac{2}{n}}= n |\Omega|^{-\frac{2}{n}},
\end{equation}
where the equality holds if and only if $L_1=\cdots=L_n$. Combining   (\ref{L1L2Ln}) and (\ref{meanL1L2}) yields that when
\begin{equation}\label{}
 |\Omega| < \left(   \frac{ nd \pi^2 }{ f'(0)+\ln(g'(0)) - \frac{|a|^2}{4d}} \right)^{n/2},
 \end{equation}
the population  must go extinct. This completes the proof.

\hfill $ \Box$

\subsection{The spectrum of the Laplacian operator} 

To find the the extreme volume size $V_{ex}$ in Theorem \ref{minVGeneral} we applied the Schwarz symmetrization argument and the classical Rayleigh-Faber-Krahn inequality and its generalizations.   In this part we discuss that one can avoid applying rearrangement type arguments to obtain inequalities for eigenvalues of the Laplacian operator with Dirichlet boundary conditions for an arbitrary domain $\Omega$.  However, the extreme volume size $V_{ex}$ that we deduce with this argument is slightly smaller than the one given in (\ref{Vex1}).   Suppose that
 \begin{eqnarray}\label{Lapeigen}
\left\{ \begin{array}{lcl}
\hfill -\Delta \phi &=&\lambda  \phi  \quad  \text{in} \ \  \Omega ,  \\
\hfill \phi&=&0   \quad  \text{on} \ \  \partial\Omega , \\
\hfill  \phi&>&0   \quad  \text{in} \ \  \Omega .
\end{array}\right.
  \end{eqnarray}
The discreteness of the spectrum of the Laplacian operator allows one to order the eigenvalues $0< \lambda_1 \le \lambda_2 \le \cdots \le \lambda_k \le \cdots $ monotonically. Proving lower bounds for $\lambda_k$ has been a celebrated problem in the field of partial differential equations.  We now  mention key results in the field and then we apply the lower bounds to our models.  In 1912, H. Weyl showed that the spectrum of (\ref{Lapeigen}) has the following asymptotic behaviour as $k\to \infty$,
\begin{equation}
\lambda_k \sim (2\pi)^2 |B_1|^{-\frac{2}{n}}  \left(\frac{k}{|\Omega|}\right)^{\frac{2}{n}}  .
\end{equation}
In 1960, P\'{o}lya in \cite{p} proved that for certain geometric shapes that is ``plane-covering domain" $\Omega \subset \mathbb R^n$  the following holds for all $k\ge 1$,
\begin{equation}
\lambda_k \ge (2\pi)^2 |B_1|^{-\frac{2}{n}}  \left(\frac{k}{|\Omega|}\right)^{\frac{2}{n}}.
\end{equation}
Then he conjectured that the above inequality should hold for general domains in $R^n$. His original proof and his conjecture were provided for $n=2$.  Even though this conjecture is still an open problem there are many interesting results in this regard.  The inequality
\begin{equation}
\lambda_k \ge  C_n (2\pi)^2 |B_1|^{-\frac{2}{n}}  \left(\frac{k}{|\Omega|}\right)^{\frac{2}{n}},
\end{equation}
with a positive constant $C_n<\frac{n}{n+2}$ was given for arbitrary domains by Lieb in \cite{l} and references therein. Li and Yau \cite{ly} improved upon this result, proving  that the following lower bound holds for all $k\ge 1$,
\begin{equation}\label{lyeig}
\lambda_k \ge  \frac{n}{n+2} (2\pi)^2 |B_1|^{-\frac{2}{n}}  \left(\frac{k}{|\Omega|}\right)^{\frac{2}{n}}.
\end{equation}

We now apply the lower bound (\ref{lyeig}) to get the following extreme volume size. Note that this $V_{ex}$ is independent from Bessel functions.
\begin{thm}\label{minVGeneralLY}
Let $A=d(\delta_{i,j})_{i,j=1}^m$ for a positive constant $d$ and the vector field $a$ be divergence free.  We assume that $f$ satisfies (F0)-(F1), $g$ satisfies (G0)-(G2) and (\ref{conditionfg}) holds.   Suppose that $\Omega\subset \mathbb R^n$ is an open bounded domain with smooth boundary and $ |\Omega| < V_{ex}$ where
\begin{equation}\label{Vex2}
V_{ex}=\Gamma\left(1+\frac{n}{2}\right) \left[\frac{4d n \pi }{ (n+2)[f'(0)+\ln(g'(0))]  }\right]^\frac{n}{2}.
\end{equation}
 Then $\lim_{m\to\infty} N_m(x)=0$ for any $x\in\Omega$.
\end{thm}

\begin{cor} In two dimensions $n=2$, the extreme volume size provided in (\ref{Vex2}) is
\begin{equation}\label{Vexpn2}
V_{ex}=\frac{2 d \pi }{f'(0)+\ln(g'(0))}. 
\end{equation}
   Moreover, in three dimensions $n=3$,  the extreme volume size (\ref{Vex2}) simplifies to
\begin{equation}\label{Vexpn3}
V_{ex}=\frac{18\sqrt 3\pi^2}{5 \sqrt 5} \left( \frac{d }{f'(0)+\ln(g'(0)) } \right)^{\frac{3}{2}} ,
\end{equation}
where we have used the fact that
$\Gamma(1+\frac{n}{2})=\Gamma(\frac{5}{2})=\frac{3\sqrt \pi}{4}$. 
\end{cor}
Note that $V_{ex}$ given in (\ref{Vexpn2}) and (\ref{Vexpn3}) are slightly smaller than the ones given by (\ref{Vexn2}) and  (\ref{Vexn3}) in Corollary \ref{coeVex2} and Corollary \ref{coeVex3}, respectively. 
\\
\\
\noindent {\it Proof of Theorem \ref{minVGeneralLY}}.  Note that we have the following lower bound on the first eigenvalue as an immediate consequence of (\ref{lyeig}),  that is
\begin{equation}\label{}
\lambda_1(dI,0, \Omega) \ge  d \frac{n}{n+2} (2\pi)^2 \left( |B_1|  |\Omega|\right)^{-\frac{2}{n}}.
\end{equation}
Following ideas provided in the proof of Theorem \ref{minVGeneral} one can complete the proof.

\hfill $ \Box$

\noindent {\bf Acknowledgement:}  The authors would like to thank the anonymous referees for many valuable suggestions.

\end{document}